\documentclass[12pt]{article}

\usepackage{amsmath,amsthm,bbm}
\usepackage{color}

\usepackage{amssymb}
\usepackage{latexsym}
\usepackage{amsfonts}

\usepackage[backrefs]{amsrefs}

\usepackage{latexsym}
\usepackage{amssymb}
\usepackage{euscript}

\def\mcc{M\raise.5ex\hbox{c}C}
\def\mccarthy{M\raise.5ex\hbox{c}Carthy}

\def\ie{{\it i.e. }}

\def\h{{\cal H}}

\def\N{{\cal N}}

\def\ei{e^{i\theta}}
\def\vare{\varepsilon}

\let\i=\infty

\def\la{\langle}
\def\ra{\rangle}
\def\={\ = \ }

\def\F{{\cal F}}

\def\C{\mathbb C}
\def\R{\mathbb R}
\def\T{\mathbb T}
\def\D{\mathbb D}

\def\be{\setcounter{equation}{\value{theorem}} \begin{equation}}
\def\ee{\end{equation} \addtocounter{theorem}{1}}
\def\beq{\begin{eqnarray*}}
\def\eeq{\end{eqnarray*}}

\def\vs{\vskip 5pt}

\def\bp{{\sc Proof: }}
\def\ep{{}{\hfill $\Box$} \vskip 5pt \par}

\def\bl{\begin{lemma}}
\def\el{\end{lemma}}
\def\bt{\begin{theorem}}
\def\et{\end{theorem}}
\def\bprop{\begin{prop}}
\def\eprop{\end{prop}}
\def\bd{\begin{definition}}
\def\ed{\end{definition}}
\def\br{\begin{remark}}
\def\er{\end{remark}}
\def\bexer{\begin{exercise}}
\def\eexer{\end{exercise}}
\def\bfig{\begin{figure}}
\def\efig{\end{figure}}

\newtheorem{theorem}{Theorem}[section]
\newtheorem{prop}[theorem]{Proposition}
\newtheorem{lemma}[theorem]{Lemma}
\newtheorem{cor}[theorem]{Corollary}

\newtheorem{conjecture}[theorem]{Conjecture}

\newtheorem{remark}[theorem]{Remark}
\newtheorem{defin}[theorem]{Definition}

\AtBeginDocument{%
 \def\MR#1{} 
}

\renewcommand\Re{\mathrm{Re\, }}
\renewcommand\Im{\mathrm{Im\, }}

\def\g{\gamma}

\def\tp{{\tilde{p}}}
\def\tq{{\tilde{q}}}
\def\tg{{\tilde{\g}}}

\newcommand\hfk{\hat{f}(k)}
\newcommand\hfn{\hat{f}(n)}
\newcommand\hg{\hat{g}}
\newcommand\hgk{\hat{g}(k)}

\newcommand\kr{Krzy\.{z}}
\newcommand{\tl}{\tau_\ell}
\newcommand{\wl}{w_\ell}
\newcommand\kn{K_n}
\newcommand\knd{K_n^\bullet}
\newcommand\ram{{\mathcal R}_m}
\newcommand\ran{{\mathcal R}_n}
\newcommand\rand{\ran^\bullet}
\newcommand\tn{{\stackrel n \sim}}
\newcommand\pn{{\mathcal P}_n}
\newcommand\pnd{{\mathcal P}_n^\bullet}

\def\pnt{\pn^{\T}}
\def\pnto{{\mathcal P}_{n;1}^{\T}}

\def\gs{\gamma^\sharp}
\def\ps{p^\sharp}

\newcommand\tlog{\big[\log | \g|^2 \big]}
\renewcommand\ei{e^{i \theta}}

\renewcommand\N{\mathbb N}
\renewcommand\O{\Omega}

\numberwithin{equation}{section}
\title{The Krzy\.{z} Conjecture and an Entropy Conjecture
}
\author{Jim Agler
\thanks{Partially supported by National Science Foundation Grant
 DMS 1665260}
\and
John E. M\raise.5ex\hbox{c}Carthy
\thanks{Partially supported by National Science Foundation Grant  
DMS 1565243
}}

\begin{document}
\maketitle

Abstract: We show that if the minimum entropy for a polynomial with roots on the unit circle is
attained by polynomials with equally spaced roots, then, under a generic hypothesis about the nature of the 
extremum, the Krzyz conjecture on the maximum modulus of the Taylor coefficients of a holomorphic
function that maps the disk to the punctured disk is true.

\setcounter{section}{0}
\section{Introduction}
Let $\O$ denote the set of holomorphic functions that map the unit disk $\D$ to
$\D\setminus \{ 0\}$.
The Krzy\.{z} conjecture, due to J.~Krzy\.{z} \cite{krz68}, is the following conjecture about the
size of Taylor coefficients of functions in $\O$.
\begin{conjecture}
\label{conk}
Let $n$ be a positive integer. Then
\be
\label{eqa1}
\knd := \sup_{f \in \O} \{ |\hfn| \ : \ f \in \O \}
\= \frac{2}{e} .
\ee
Moreover, equality is obtained in \eqref{eqa1} only for functions of the form
\be
\label{eqa2}
f(z) \= \zeta \exp( \frac{ z^n + \omega}{ z^n - \omega })
\ee
where $\zeta$ and $\omega$ are unimodular constants.
\end{conjecture}
For any function $f$ defined and holomorphic on a neighborhood of the origin,
we use $\hfk$ to denote the $k^{\rm th}$ Taylor coefficient at $0$, so
\[
\hfk \= \frac{f^{(k)}(0)}{k!} .
\]
For a history of the Krzy\.{z} conjecture and a summary of known results, see Section~\ref{secb}.

The purpose of this note is to establish a connection between the Krzy\.{z} conjecture and the following conjecture
about the entropy of polynomials with roots on the unit circle $\T$.
\begin{conjecture}
\label{cone}
Let $p$ be a non-constant polynomial, all of whose roots lie on $\T$, and normalized so that
$\frac{1}{2 \pi} \int_0^{2\pi} |p(e^{i\theta} )|^2 d\theta = 1$.
Then
\be
\label{eqa3}
\frac{1}{2 \pi} \int_0^{2\pi}  |p(e^{i\theta} )|^2 \log |p(e^{i\theta} )|^2 d\theta
\ \geq \ 1 - \log(2).
\ee
Moreover, equality occurs in \eqref{eqa3} only for polynomials of the form
\be
\label{eqa4}
p(z) \= \frac{\zeta}{\sqrt{2}} ( \omega + z^n) ,
\ee
where $\zeta$ and $\omega$ are unimodular constants, and $n$ is a positive integer.
\end{conjecture}

\vs
We shall let $\h$ denote the Herglotz class, the holomorphic functions
on the unit disk that have non-negative real part.
A function $f$ is in $\Omega$ if and only if there is a function $g \in \h$ so that
\[
f \= e^{-g} .
\]
Given an $(n+1)$-tuple $a = (a_0, \dots, a_n)$ of complex numbers, we shall say that $a$ is {\em solvable Herglotz data}
if there exists $g \in \h$ satisfying
\[
\hgk \= a_k, \ {\rm for\ } k = 0, \dots, n .
\]
We shall say that $a$ is {\em extremal Herglotz data} if
it is solvable but for any $r > 1$, the data $(a_0, r a_1, \dots, r^n a_n)$ is not solvable.

We shall say that $f$ is $\knd$-extremal if $f$ is in $\Omega$ and $\hfn = \knd$.
It was proved in \cite{hsz}
 that 
if $f$ is $\knd$-extremal, and $g = - \log(f)$,
then the first $n+1$ Taylor coefficients of $g$ are extremal Herglotz data (we give
a proof of this in Lemma~\ref{lemc2}).
By a theorem of G. Pick \cite{pi16}, this means $g$ must have the form
\be
\label{eqa5}
g(z) \= ai + \sum_{\ell =1}^m \wl \frac{\tl + z}{\tl -z} ,
\ee
where $a \in \R$, the number $m$ satisfies $ 1 \leq m \leq n$, each $\wl > 0$,
and each $\tl$ is a distinct point on $\T$.
We shall let $\ram$ denote the set of rational functions that have the form \eqref{eqa5}
(that is, rational functions of degree $m$ that are $m$-fold covers of the right-half plane by the unit disk),
and we shall let $\rand = \cup_{1 \leq m \leq n} \ram$.

Our first main result analyzes the critical points for the \kr\ functional.
We prove in Theorem~\ref{thme1} that if
$g$ is in $\ran$ (and is normalized in a way described in Section~\ref{secf}), then
 $f = e^{-g}$ is a critical point for $\knd$ with critical value $\eta$
 if and only if 
\[
e^{-g}\  \tn\  \eta \g^2 ,
\]
where $\gamma$ is a polynomial of degree $n$ of unit norm in the kernel of $\Re g(S_n)$,
where $S_n$ is the compression of the unilateral shift to polynomials of degree less than or equal to $n$, and
the notation $\tn$ means that the functions have the same Taylor coefficients up to degree $n$.

We use this result to prove  that if the extremals $f$  for the \kr\ problem have $g= -\log f$ of full degree,
then the entropy conjecture implies the \kr\ conjecture.
\bt
\label{thma1}
Let $f$ be $\knd$-extremal, and assume that $g = - \log f$ is in $\ran$. If Conjecture~\ref{cone} is true,
then $f$ has the form \eqref{eqa2}.
\et
We prove Theorem~\ref{thma1} in Section~\ref{secg}.
In Section~\ref{sech} we study critical points of the entropy functional from \eqref{eqa3}.
In Section~\ref{seci}, we prove a special case of Conjecture \ref{cone}.
Finally, in Section~\ref{secj}, we show how Conjecture \ref{cone} would follow
from Conjecture \ref{conb}, due to A. Baernstein II.

\section{History of the Krzy\.{z} conjecture}
\label{secb}
J. \kr\ proved Conjecture \ref{conk} for $n=2$, and conjectured it for all $n$.
The $n=3$ case was proved by J. Hummel, S. Scheinberg and L. Zalcman \cite{hsz};
they also proved Lemma~\ref{lemc2} below, and that \eqref{eqa2} is a strict local minimum for
\eqref{eqa1} (after normalizing so that $f(0)$ and $\hat{f}(n)$ are both positive).
The $n=4$ case was proved first by D. Tan \cite{ta83}, and later by a different method by J. Brown  \cite{br87}.
The $n=5$ case was proved by N. Samaris \cite{sam03}.

C.~Horowitz \cite{hor78} proved that there is some constant $H < 1$ such that
$\knd \leq H$ for all $n$; his proof showed $H \leq  0.99987\dots$. This was improved by 
R.~Ermers to $H \leq 0.9991\dots$ \cite{erm90}.

In \cite{msuv15}, M. Martin,  E. Sawyer, I. Uriarte-Tuero,
and D. Vukoti\'c
prove that 16 different conditions are all equivalent to the \kr\ conjecture.
The paper also includes a useful historical summary.

\vs
Conjecture~\ref{cone} may be compared with the following sharp inequality, conjectured by I. Hirschman \cite{hi57}
and proved by W. Beckner \cite{be75}:

If $f \in L^2(\R)$ has norm $1$, and $\F f$ denotes the Fourier transform of $f$, then
\be
\label{eqb1}
\int |f|^2 \log |f|^2 \ + \ \int |\F f|^2 \log |\F f|^2 \ \leq \ \log(2) -1 .
\ee
Equality is obtained in \eqref{eqb1} for Gaussians.

A. Baernstein II made the following conjecture in 2008 \cite{ba08},
where the quasi-norms are with respect to normalized Lebesgue measure on the circle.
By $\| f \|_0$ we mean $\exp ( \int_\T \log |f| )$.

\begin{conjecture}
\label{conb}
Let $Q(z) = 1 + z^n$. Then for all $0 \leq s \leq t \leq \infty$, and for all non-constant polynomials
$p$ with all their roots on the unit circle, 
\be
\label{eqb2}
\frac{\| p \|_s}{\| Q \|_s} \ \leq \ \frac{\| p \|_t}{\| Q \|_t}.
\ee
\end{conjecture}

In Section~\ref{secj} we show how Baernstein's conjecture implies the entropy conjecture.

\section{Preliminaries}
\label{secc}
Suppose $f$ and $g$ are analytic functions on a neighborhood of $0$, and $n \in \N$. Say
\[
f \tn g 
\]
if $\hfk = \hgk$ for $0 \leq k \leq n$.
We leave the proof of the following lemma to the reader.
\bl
\label{lemc1}
Suppose $f$ and $g$ are analytic on a neighborhood of zero. Assume that $f(0) = g(0) = b$,
and $\phi$ is analytic in a neighborhood of $b$. If $f \tn g$, then $\phi \circ f \tn \phi \circ g$.
\el
The following result is due to G. Pick \cite{pi16}, and can be found in any book on Pick interpolation
such as \cite{bgr90, foi-fra, ampi}.
\bl
\label{lemc3}
If $a = (a_0, \dots , a_n)$ is extremal Herglotz data, then there exists a unique function $g \in \h$ such 
that $\hgk = a_k,$ for $0 \leq k \leq n$. Moreover, $g \in \rand$. Conversely, if $g \in \rand$, then
$a = (\hg(0), \dots, \hg(n))$ is extremal Herglotz data.
\el
\bl
\label{lemc2}
Fix $n \geq 1$ and assume that $f$ is $\knd$-extremal. Define $g$ by $f = e^{-g}$.
Then $g \in \rand$.
\el
\bp
We know that $g$ must be in $\h$, so by Lemma~\ref{lemc3}, if $g$ is not in $\rand$, then for
some $r > 1$ we have a function $h \in \h$ such that 
\[
\hat{h}(k) \=  r^k \hgk, \quad 0 \leq k \leq n .
\]
Then $\phi = e^{-h}$ is in $\O$, and by Lemma~\ref{lemc1},
$\phi(z) \tn f(rz)$, so
\[
| \hat{\phi}(n) | \= r^n | \hfn| \ > \ |\hfn|.
\]
This contradicts the claim that $f$ is extremal.
\ep
\begin{cor}
\label{corc1}
\[
\knd \= \sup_{g \in \rand} |\widehat{(e^{-g})} \ (n)|.
\]
\end{cor}
Given that $\rand = \cup_{m=0}^n \ram$, Corollary~\ref{corc1} suggests the following optimization problem.
For each $n \geq 1$, define
\[
\kn \= \sup_{g \in \ran} |\widehat{(e^{-g})} \ (n)|. 
\]
\begin{remark}
{\rm
As $\ran$ is a dense open set in $\rand$, we have $\kn = \knd$. However, whereas a normal families argument guarantees that an extremal function for $\knd$ always exists, it is not obvious that an extremal for $\kn$ exists.
If \kr's conjecture is true, then the supremum is attained.
}
\end{remark}

\section{The critical points of $\kn$}
\label{secd}

\subsection{The definition of critical points}
For the rest of the paper, $n$ will be a positive integer.
There are a number of equivalent ways to view $\ran$ as a topological vector space.
\begin{enumerate}
\item
Using the $a, w, \tau$ parameters of \eqref{eqa5} (where $m=n$).
\item
As the subset of the space of extremal Herglotz data points $a = (a_0, \dots, a_n)$ with the property
that $(a_0, \dots, a_{n-1})$ is not extremal.
\item
With the topology of uniform convergence on compact subsets of $\D$.
\end{enumerate}
We would like to consider the local maxima of the function $F : \ran \to \R$ defined by
\be
\label{eqd2}
F(g) \= | \widehat{(e^{-g})}(n) |^2 .
\ee
Let  $\pnd$ denote the set of complex polynomials of degree less than or equal to $n$, and $\pn$ the polynomials
of degree exactly $n$.
\begin{defin}{\rm
Let $g \in \ran$. We say $d$ is an {\em admissible direction} at $g$ if $d \in \pnd$ and there exists $\vare >0$
 such that
 \[
 (\hat{g}(0) + t \hat{d}(0), \dots, \hat{g}(n) + t \hat{d}(n))
 \]
 is solvable Herglotz data for all $t$ in $(0,\vare)$.
We say that $g$ is a critical point for $\kn$ if 
  \be
  \label{eqd1}
    \frac{d}{dt} |  \widehat{e^{ -(g + t d)}}(n)  |^2 \big|_{t = 0^+} \le 0
\ee
  whenever $d$ is an admissible direction at $g$. 
 If $g$ is a critical point for $\kn$, then we refer to $\eta = \widehat{e^{-g}}(n)$ as the {\em critical value}.
}
  \end{defin}

\subsection{A Hilbert space setting for the analysis of critical points}

Let $H^2$ denote the classical Hardy space on the unit disk. We shall think of $\pnd$ as a subspace
of $H^2$,  and let $P_n$ be the orthogonal projection
from $H^2$ onto $\pnd$. Define an operator $S_n$ on $\pnd$ by
the formula 
\[
(S_n q )(z)  \= P_n (z q(z)), \qquad q(z) \in \pnd .
\]
The operator $S_n$ is the truncated shift, and is nilpotent of order $n+1$. Hence if $f$ is any holomorphic
function on a neighborhood of $0$, we can define $f(S_n)$ by the Riesz functional calculus, or by either of the 
two equivalent formulas
\beq
f(S_n) &\=& \sum_{k=0}^n \hfk s_n^k \\
f(S_n) q &=& P_n( f q) .
\eeq
Observe that if $f$ and $g$ are both holomorphic on a neighborhood of $0$, then 
\[
f \tn g \ \Leftrightarrow \ f(S_n) = g(S_n).
\]
The following two propositions are basically a reformulation of Lemma~\ref{lemc3}
to the Hilbert space interpretation of interpolation.
\begin{prop}
\label{propd1}
Let $a$ be an $(n+1)$-tuple of complex numbers.
Then $a$ is solvable Herglotz data if and only if
\[
\Re \sum_{k=0}^n a_k S_n^k \geq 0 .
\]
Moreover $a$ is extremal Herglotz data if and only if $\Re \sum_{k=0}^n a_k S_n^k $ is positive semi-definite
but not positive definite.
\end{prop}
\begin{prop}
\label{propd2}
Assume $g \in \h$.
\begin{enumerate}
\item $\Re g(S_n) \geq 0 .$
\item The function $g$ is in $\rand$ if and only if $\Re g(S_n) $ is singular.
\item If $0 \leq m \leq n$, then $g \in \ram$ if and only if ${\rm rank} (g(S_n) ) = m$.
\end{enumerate}
\end{prop}

\subsection{Local maxima are critical points}

We need to show that there is enough smoothness at local maxima to make sense of \eqref{eqd1},
at least when the local maximum is in $\ran$.
\begin{prop}
\label{propd3}
Let $F$ be defined by \eqref{eqd2}. If $g \in \ran$ is a local maximum for $F$, then $g$
is a critical point for $\kn$.
\end{prop}
\bp
Let $d$ be an admissible direction for $g$. Thus by Proposition~\ref{propd1} there exists $\vare > 0$ so that
\[
\Re [ g(S_n) + t d(S_n) ] \ \geq \ 0 \quad \forall\  t \in [0, 
\vare) .
\]
Let $\rho(t)$ denote the smallest eigenvalue of $\Re[g(S_n) + td(S_n)]$, so for each
$t$ we have $\Re[g(S_n) + td(S_n) - \rho(t)]$ is positive semidefinite and singular.
By Proposition~\ref{propd1}, for every $t$  there exists $g_t \in \rand$ such that 
\[
g_t\  \tn \ [g + td - \rho(t) ].
\]
Since $g$ is a local maximum, we have $\rho(t) \to 0$ as $t \to 0^+$, and
$g_t \to g$. As $\ran$ is open in $\rand$, this means for some 
 $\delta > 0$, we have $g_t \in \ran$ for all $t$ in $[0,\delta)$.
 As $g$ is a local maximum for $F$, we have
 \[
 F(g_t) \ \leq \ F(g) \quad \forall \ t \in [0,\delta) .
\]
As
\[
F(g+td) \= | \widehat{e^{-(g+td)} }(n) |^2 ,
\]
it is differentiable with respect to $t$, and as
\[
F(g+td) \= e^{-2 \rho(t)} F(g_t) \ \leq \ F(g) ,
\]
this derivative is nonpositive at $0$.
\ep

\subsection{Some lemmas about critical points}
\label{secd}

We shall let $\| \gamma \|$ denote the $H^2$-norm, so
\[
\| \gamma \|^2 \= \frac{1}{2\pi} \int_{0}^{2\pi} | \gamma(e^{i\theta}) |^2 d \theta .
\]
We shall let $n$ be fixed, and write $S$ for $S_n$ for legibility.
\bl
\label{lemd1}
If $g \in \ran$, then there exists a unique vector $\gamma$ in $\pn$ such that $\hat{\gamma}(n) > 0$, 
$\| \gamma\| =1$ and $\ \Re g(S) \gamma = 0$. Furthermore, if
\[
g(z) \= ai + \sum_{\ell =1}^n \wl \frac{\tl + z}{\tl -z} ,
\]
then
\[
\gamma(z) \= \frac{1}{\nu} \prod_{\ell =1}^n (z - \tau_l) ,
\]
where
\[
\nu \= \|  \prod_{\ell =1}^n (z - \tau_l) \|.
\]
\el
\bp
By Proposition~\ref{propd2}, ${\rm rank}(g(S) = n$. Hence there exists some nonzero vector $q$ in ${\rm ker\ }\Re g(S)$.
We will show that $q(\tau_\ell) = 0$ for each $\ell$, and then define
\[
\gamma(z) \= \frac{|\hat{q}(n)|}{\hat{q}(n) \| q \|} q .
\]
As
\beq
\Re \frac{\tau + S}{\tau - S} &\=&
\frac{1}{2} \left( \frac{\tau + S}{\tau - S} + \frac{\bar \tau + S^*}{\bar \tau - S^*} \right) \\
&=&
(\bar \tau - S^*)^{-1} (1 - S^* S) (\tau - S)^{-1} \\
&=&
(\bar \tau - S^*)^{-1} (z^n \otimes z^n) (\tau - S)^{-1} \\
&=&
[(\bar \tau - S^*)^{-1} z^n] \otimes [ (\bar \tau - S^*)^{-1} z^n],
\eeq
we have
\[
\Re g(S) \= \sum_{\ell =1}^n \wl [(\bar \tl - S^*)^{-1} z^n] \otimes [ (\bar \tl - S^*)^{-1} z^n].
\]
Since each $\wl >0$, we can only have $\la \Re g(S) q, q \ra = 0$ if for each $\ell =1 , \dots, n$ we have
\[
\la q,  (\bar \tl - S^*)^{-1} z^n \ra = 0.
\]
As 
\[
 (\bar \tl - S^*)^{-1} z^n \= \tau^{n+1}( 1 +  \bar \tau z + \dots + \bar \tau^n z^n ),
 \]
 we get
\beq
 \la q,  (\bar \tl - S^*)^{-1} z^n \ra &\= &
 \bar \tau^{n+1} \la \hat{q}(0) + \dots + \hat{q}(n)z^n, 1 + \bar \tau z + \dots \bar \tau^n z^n \ra \\
 &=&
 \bar \tau^{n+1} q(\tau) .
 \eeq
Therefore $\la \Re g(S) q, q \ra = 0$ implies that $q$ vanishes at each $\tau_\ell$, as claimed.
\ep
\bl
\label{lemd2}
Let $g \in \ran$ and let $\gamma$ be the vector described in Lemma~\ref{lemd1}.
For $d \in \pnd$, the following hold:
\begin{enumerate}
\item If $d$ is an admissible direction at $g$, then $\la \Re d(s) \g, \g \ra \geq 0$.
\item
If $\la \Re d(s) \g, \g \ra = 0$, then $d$ is an admissible direction at $g$ if and only if $\ \Re d(s) \g = 0$.
\item
If $\la \Re d(s) \g, \g \ra = 0$, then $d + \vare$ is an admissible direction for every $\vare > 0$.
\end{enumerate}
\el
\bp
1. Since $g + td \in \h$ for $t$ small and positive, we must have 
\[
\la \Re[ g(S) + t d(S)]\g, \g \ra \= t \la \Re d(S) \g, \g \ra \ \geq \ 0 .
\]
2. If $\beta \perp \gamma$, then 
\[
\la  \Re[ g(S) + t d(S)] (\g + \beta) , ( \g + \beta) \ra
\=2 t \la \Re d(S) \gamma , \beta \ra  .
\]
The right-hand side is non-negative for all $\beta$ and small positive $t$ if and only if
$\Re d(S) \gamma = 0$.

3. If $\beta \perp \gamma$, then 
\[
\la  \Re[ g(S) + t( d(S)+ \vare)] (\g + \beta) , ( \g + \beta) \ra
\=2 t \la \Re d(S) \gamma , \beta \ra   + t \vare (\| \g \|^2 + \| \beta \|^2).
\]
By Cauchy-Schwarz, for any fixed $\vare > 0$, the right-hand side is non-negative for $t$ positive and sufficiently small.
\ep
\bl
\label{lemd3}
Let $g \in \ran$ and let $\gamma$ be the vector described in Lemma~\ref{lemd1}.
If $g$ is a critical point for $\kn$ and $d \in \pnd$ satisfies $\Re \la d(S) \g, \g \ra = 0$,
then
\[
\frac{d}{dt} F(g + td) |_{t=0^+} \ \leq \ 0.
\]
\el
\bp
By Lemma~\ref{lemd2}, for all $\vare > 0$ we have
$d + \vare $ is admissible, so by Proposition~\ref{propd3} we have
\[
\frac{d}{dt} F(g + t(d+ \vare)) |_{t=0^+} \ \leq \ 0.
\]
Now let $\vare \to 0^+$. 
\ep
\bl
\label{lemd4}
If $g$ is analytic on a neighborhood of $0$ and $d \in \pnd$, then
\[
\frac{d}{dt} | \hat{e^{-g+td}}(n) |^2 \Big|_{t=0} \= -2 \Re \la z^n , e^{-g(S)} 1 \ra \la d(S) e^{-g(S)} 1, z^n \ra .
\]
\el
\bp
Computation.
\ep

\section{The Critical Point Equation}
\label{sece}

We fix $n \geq 1$, and write $S$ for $S_n$. For $p \in \pnd$, define $\tp$ by 
\[
\tp(z) = z^n \overline{p(\frac{1}{\bar z})} .
\]
So if $p(z) = a_0 + a_1 z+ \dots + a_n z^n$, then 
$\tp(z) = \bar a_n + \bar a_{n-1} z + \dots + \bar a_0 z^n$.
We shall say $p$ is {\em self-inversive} if $p = \tp$.

\bt
\label{thme1}
Let $g \in \ran$ and let $\gamma$ be the vector described in Lemma~\ref{lemd1}.
Then  $g$ is a critical point of $\kn$ with critical value $\eta$ if and only if
\be
\label{eqe2}
e^{-g}\  \tn\  \eta \g \tg .
\ee
\et
\bp
Suppose $g$ is a critical point of $\kn$ with critical value $\eta$.
If  $d \in \pnd$ and $\Re \la d(S) \g, \g \ra = 0$, then by Lemma~\ref{lemd3},
\[
\frac{d}{dt} | \la e^{-(g+td)}(S) 1, z^n \ra |^2 \Big|_{t=0^+} \= \frac{d}{dt} F(g+td) \Big|_{t=0^+} \ \leq \ 0.
\]
Hence if $\la d(S) \g , \g \ra = 0$, so 
$\Re \la \zeta d(S) \g, \g \ra = 0$ for all $\zeta \in \T$, we get by Lemma~\ref{lemd4}
\[
- 2 \Re \zeta \bar \eta \la d(S) e^{-g}(S) 1, z^n \ra \ \leq \ 0 .
\]
As this holds for all $\zeta \in \T$, we get that
\be
\label{eqe3}
d \in \pnd \ {\rm and\ } \la d(S) \g , \g \ra = 0 \
\Rightarrow \ \la d(S) e^{-g}(S) 1, z^n \ra = 0 .
\ee
Equivalently,
\[
\sum_{k=0}^n \hat{d}(k) \la S^k \g , \g \ra = 0 \ \Rightarrow \
\sum_{k=0}^n \hat{d}(k) \la S^k e^{-g}(S) 1 , z^n \ra = 0 .
\]
By duality (in the finite dimensional space $\pnd$), this means there exists $c \in \C$ so that 
\be
\label{eqe1}
  \la S^k e^{-g}(S) 1 , z^n \ra \= c \ \la S^k \g , \g \ra , \quad 0 \leq k \leq n.
  \ee
Letting $k=0$ in \eqref{eqe1}, we get $c = \eta$.
So for $0 \leq k \leq n$, \eqref{eqe1} gives 
\beq
\widehat{e^{-g}}(n-k) &\=&\la S^k e^{-g}(S) 1 , z^n \ra  \\
&=& \eta \la S^k \g , \g \ra\\
&=&
\eta \la z^k \g, \g \ra_{H^2} \\
&=& \eta \int_{0}^{2\pi} e^{ik\theta} \g(e^{i\theta}) \overline{ \g(e^{i\theta}) } \frac{d\theta}{2\pi} \\
&=& \eta \int_{0}^{2\pi} e^{ik\theta} \g(e^{i\theta}) [e^{-in\theta}  \tg(e^{i\theta}) ] \frac{d\theta}{2\pi} \\
&=& \eta \widehat{\g \tg}(n-k).
\eeq
Therefore $e^{-g}\  \tn\  \eta \g \tg $, as desired.

Conversely, suppose \eqref{eqe2} holds. Reversing the logic, we get that 
\eqref{eqe3} holds. This means that on the $n$-dimensional subspace of $\pnd$
given by
\be
\label{eqe4}
\{ q : \la q(S) \g , \g \ra = 0 \} ,
\ee
we have
\[
\frac{d}{dt} F(g + td) \Big|_{t=0} \= 0 .
\]
The orthocomplement of the subspace \eqref{eqe4} is spanned by the function
$\beta = P_n ( |\g|^2)$.
An arbitrary polynomial in $\pnd$ can be written as $q + a \beta$,
where $q$ is in \eqref{eqe4} and $a \in \C$.
By Lemma~\ref{lemd2}, this is an admissible direction if and only if $\Re a \geq 0$.
By Lemma~\ref{lemd4}, and using \eqref{eqe2}, 
\beq
\frac{d}{dt} F( g + ta \beta) \Big|_{t=0^+} &\=&
- 2 \Re a  \la z^n , \eta \g \tg \ra \la \beta \eta \g \tg , z^n \ra  .
\\
&=&
- 2 \Re a |\eta|^2 |\la z^n , \g \tg \ra |^2
\eeq
and this is less than or equal to $0$ whenever  $\Re a \geq 0$.
This means $g$ is a $\kn$ critical point.
Finally, 
\beq
\widehat{e^{-g}}(n) &\= &  \la \eta \g \tg , z^n \ra \\
&=&
\eta \la \g z^n \bar \g , z^n \ra \\
&=& \eta .
\eeq
\ep

\section{Normalization}
\label{secf}

If $g$ is as in \eqref{eqa5} and is a local maximum for $F$, then $bi + g(\zeta z)$ is also a local maximum
for any unimodular $\zeta$ and real $b$. We can choose $\zeta$ and $b$ so 
that $a=0$ and $\prod_\ell (- \tau_\ell) = 1$.
\begin{defin}
If $g \in \ran$, we say that $g$ is {\em normalized} if $g$ has the form
\be
\label{eqf1}
g(z) \=  \sum_{\ell =1}^n \wl \frac{\tl + z}{\tl -z} ,
\ee
where $\prod_{\ell =1}^n  (- \tau_\ell) = 1$.
\end{defin}
\bl
\label{lemf1}
Let $g \in \ran$ be a local maximum for $F$ with critical value $\eta$ and let $\gamma$ be as in 
 Lemma~\ref{lemd1}. Assume $\Im g(0)$ is chosen in the range $[-\pi,\pi)$.
 Then $g$ is normalized if and only if $\eta > 0$ and
 $\g $ is self-inversive.
 \el
 \bp
 We have
 \[
 \prod_\ell (z- \tl)  \= \prod_{\ell =1}^n  (- \tau_\ell)  \ [ \prod_\ell (z- \tl) ]\ \tilde{ } .
 \]
 So by Lemma~\ref{lemd1}, $\g = \tg$ if and only if $ \prod_\ell (-\tl) = 1$.
From \eqref{eqa5}, we have
\[
g(0) \= ai + \sum_{\ell =1}^n w_\ell ,
\]
and from Theorem~\ref{thme1}, we have 
\[
e^{-g(0)} \= \eta \g(0) \tg(0) \= \frac{\eta}{\nu^2}   \prod_\ell (-\tl) .
\]
So if $ \prod_\ell (-\tl) = 1$, then $\eta$ is positive if and only if $\Im g(0)$ is a multiple of $2\pi$.
\ep
\begin{prop}
\label{propf1}
If $g$ is a normalized local maximum for $F$ with critical value $\eta$ then
\be
\label{eqf2}
-g \ \tn \ \log\frac{\eta}{\nu^2} + \log \left( \prod_{\ell=1}^n (1 - \bar \tl z)^2 \right) .
\ee
\end{prop}
\bp
Since $g$ is normalized, we have
\[
\g(z) \tg(z) \=
\frac{1}{\nu^2} \prod_{\ell =1}^n (z - \tl) \prod_{\ell =1}^n (1 - \bar \tl z)  \=
 \frac{1}{\nu^2}\prod_{\ell =1}^n (1 - \bar \tl z)^2 .
 \]
 So from Theorem~\ref{thme1}, we have
 \[
 e^{-g} \ \tn \  \frac{\eta}{\nu^2}\prod_{\ell =1}^n (1 - \bar \tl z)^2 .
 \]
 Then \eqref{eqf2} follows from Lemma~\ref{lemc1}.
 \ep
 \begin{prop}
 \label{propf2}
 If $g$ is a normalized local maximum for $F$ with critical value $\eta$ then
 \[
 \sum_{\ell=1}^n w_\ell \= \log \frac{\eta}{\nu^2} .
 \]
 For $k = 1, \dots, n$,
 \[
  \sum_{\ell=1}^n w_\ell \tl^k \= \frac{1}{k}  \sum_{\ell=1}^n \tl^k.
  \]
  \end{prop}
\bp
Expand both sides of \eqref{eqf2} into power series and equate coefficients.
From \eqref{eqf1},
\[
-g(z) \=
- \sum_{\ell=1}^n w_\ell - 2 \sum_{k=1^\i} \left( \sum_{\ell =1}^n w_\ell \bar \tl^k \right) z^k .
\]

We have
\[
 \log \left( \prod_{\ell=1}^n (1 - \bar \tl z)^2 \right) 
 \=
 - 2 \sum_{k=1^\i} \left( \sum_{\ell =1}^n \frac{ \bar \tl^k }{k}\right) z^k .
 \]
 Comparing these we get the result.
 \ep
 
\section{The proof of Theorem~\ref{thma1}}
\label{secg}
Let us assume that the Entropy Conjecture~\ref{cone} holds, and that $g \in \ran$ is a normalized local 
maximum for $F$ with critical value $\eta$. Let $\gamma$ be as in Lemma~\ref{lemd1}.

By Proposition~\ref{propd3}, $g$ is a critical point for $\kn$, and by Theorem~\ref{thme1} and Proposition~\ref{propf1},
 the $(n+1)$-by-$(n+1)$ Toeplitz matrix $-\Re g(S_n)$ is the same as
the Toeplitz matrix on $\pnd$ whose entries come from the Fourier series of $\log \eta + \log |\g|^2$.
In particular, for any polynomials $p,q \in \pnd$, we have
\be
\label{eqg1}
\la - \Re g(S_n) p, q \ra \=
\frac{1}{2\pi}  \int_0^{2\pi} [ \log \eta + \log | \g(e^{i\theta}) |^2] p(e^{i\theta}) \overline{q(e^{i\theta})} d\theta.
\ee
Let $p=q =\g$, and observe that the left-hand side of \eqref{eqg1} vanishes, so 
\[
\frac{1}{2\pi}  \int_0^{2\pi} | \g(e^{i\theta})|^2  \log | \g(e^{i\theta}) |^2] d\theta
\= - \log \eta .
\]
If Conjecture~\ref{cone} holds, then
\[
- \log \eta \geq 1 - \log 2 ,
\]
and so $\eta \leq \frac{2}{e}$.
\ep

\section{Entropy Conjecture}
\label{sech}

Let us establish some notation. We shall fix $n \geq 1$ a positive integer.
All integrals are integrals over the unit circle with respect to normalized Lebesgue measure, 
and norms and inner products  are in $L^2(\T)$  with respect to this measure.

If $p$ is self-inversive, then all its zeroes either occur on $\T$, the unit circle, 
or occur in pairs $(\zeta, 1/\bar \zeta)$, or occur at the origin if ${\rm deg}(p) < n$.
We shall let $\pnt$ denote the set of polynomials in $\pn$ that are self-inversive and have all their zeroes on $\T$, and $\pnto$ denote the unit sphere of $\pnt$, {\em viz.} the polynomials
of norm $1$ in $\pnt$. 

We shall let $\Pi_n$ be the orthogonal projection 
from $L^2(\T)$ onto $\pnd$, \ie
\[
\Pi_n (\sum_{-\i}^\i c_k z^k ) \=
\sum_{0}^n c_k z^k.
\]
If $f, g$ are in $L^2(\T)$, we shall write $f \sim g $ to mean $\Pi_n f = \Pi_n g$, \ie
\[
f \sim g \ \Leftrightarrow \ \hat f(k) = \hat g(k) \ \forall \ 0 \leq k \leq n.
\]

Note that if $p,q \in \pnd$, then 
\[
\la \tp , \tq \ra \= \la q , p \ra ;
\]
in particular if $p$ and $q$ are both self-inversive, then their inner product is real.

If $p \in \pnd$, then $\ps$ denotes the polynomial
\[
\ps(z) = i \left( p(z) - \frac{2}{n} z p'(z) \right).
\]
 In terms of Fourier coefficients,
\[
\hat \ps(k) \= i  \left(\frac{n-2k}{n} \right) \hat p(k) .
\]
\begin{lemma}
\label{lemhh1}
If $p$ is in $\pnt$, then so is $\ps$.
\end{lemma}
\bp
A calculation shows that for $p$ self-inversive
\be
\label{eqhh1}
\frac{d}{d\theta} |p(e^{i\theta})|^2 \= -n \bar z^n p (z) \ps(z) .
\ee
So  $\ps$  has zeroes at the local maxima
of $|p|^2$ on $\T$; these interleave the zeroes of $p$.
If $p$ has a zero of order $k > 1$ at some point $\tau$ on $\T$, then 
\eqref{eqhh1} vanishes to order $2k-1$ at $\tau$, so $\ps$ has a zero of order $k-1$.
Counting them all up, we get that $\ps$ has $n$ zeroes on $\T$, and since it is of degree $n$,
it must be in $\pnt$.
\ep

Observe that
$\la p , \ps \ra = 0$ for all self-inversive  $p$ in  $\pn$.

Let \[
\F(p) \= \int |p|^2 \log |p|^2 .
\]
If $\gamma$ is norm one and $\F(\gamma) = m$, then 
minimizing $\F(c \gamma)$ over all $c \geq 0$,  one gets that
\be
\label{eqh18}
\F( c \gamma) \ \geq \ -  e^{-1-m} ,
\ee
with equality when $c^2 = e^{-1-m}$.
The entropy conjecture \ref{cone} is equivalent to the conjecture that the minimum of $\F(p)$ over
all $p \in \pnt$ (not just those of norm one) is $-2 e^{-2}$,
and that, up to the normalization of requiring that $\hat{p}(0)$
and $\hat{p}(n)$ are positive, this value is attained 
uniquely by the polynomial
\be
\label{eqhh5}
p(z) \=  e^{-1} (1 + z^n) .
\ee

For any function $f$ in $L^2$, we shall let $\big[ f \big]$ denote the $(n+1)$-by-$(n+1)$   Toeplitz matrix
with $(i,j)$ entry $\hat{f}(j-i)$. We shall think of this as acting on $\pnd$. If $p(z) = c \prod_{\ell=1}^n ( 1 - \bar \tl z)$, then
\[
\big[ \log |p|^2 \big] \=
\begin{bmatrix}
\log |c|^2 & - \sum \tl &  - \frac{1}{2} \sum \tl^2 \ &
\cdots & - \frac{1}{n} \sum \tl^n\\
\\
- \sum \bar \tl & \log |c|^2 & - \sum \tl & \cdots & - \frac{1}{n-1} \sum \tl^{n-1}\\
\vdots & \vdots & \vdots & \vdots & \vdots \\
- \frac{1}{n} \sum \bar\tl^n &  - \frac{1}{n-1} \sum \bar \tl^{n-1} &  - \frac{1}{n-2} \sum \bar \tl^{n-2} &
\cdots & \log |c|^2
\end{bmatrix}
\]

\bt
\label{thmh1}
Suppose $\g$ is a local minimum for $\F$ on $\pnt$, and that all the zeroes of $\g$ are distinct.
Then
\be
\label{eqhh2}
\tlog \g \= - \g .
\ee
\be
\label{eqhh3}
\tlog \ \geq - 3 .
\ee
\et
\bp
Since all the zeroes of $\g$ are distinct, it $q$ is any self-inversive polynomial in $\pnd$, then for $t$ small and real,
$\g + tq$ is self-inversive, and the zeroes must be close to the zeroes of $\g$, so they must all lie on the circle.
Therefore if we expand $\F(\g + tq)$ in powers of $t$, the first order term must vanish, since $\g$ is a critical point,
and the coefficient of $t^2$ must be non-negative, since $\g$ is a local minimum.

Calculating,  using the fact that if $p$ is self-inversive, then on the unit circle  $ z^n \overline{p(z)} = p(z)$, 
and writing 
\[ \log | \g + tq|^2 = \log |\g|^2 \ +\   2 \Re t \frac{q}{\g}\  - \ \Re t^2 \frac{q^2}{\g^2} \ + \ O(t^3), \]
 we get,
 \beq
 \F(\g + tq) &=& 
\int \big(  \log |\g|^2  +   2 \Re t \frac{q}{\g}  -  \Re t^2 \frac{q^2}{\g^2} \big)
\big( |\g|^2 + 
2 \Re t \g \bar q + t^2 |q|^2 \big) \\
&&\quad \ +\  O(t^3)   \\
&=& \F(\g) \ +\  t  \Big( 2\Re  \la \tlog \g, q \ra + 2 \Re \la \g, q \ra \Big)\\
 &&\quad  + \ t^2 \Big( \la \tlog q , q \ra + 4  \la q, q \ra - \la q, q \ra \Big)\   \ +\  O(t^3) .
 \eeq

Since at a critical point the coefficient of $t$ must vanish for all $q$, we get
$\tlog \g + \g = 0$, giving \eqref{eqhh2}. The non-negativity of the coefficient of $t^2$ gives \eqref{eqhh3}.
\ep
At a critical point, $\tlog$ will have one eigenvalue equal to $-3$, so the inequality in \eqref{eqhh3} cannot be strict.
\bprop
\label{proph1}
Suppose $\g $ is in $\pnt$ and 
\be
\label{eqhh4}
\tlog \g = \kappa \g .
\ee 
Then
\be
\label{eqhh6}
\tlog \gs \= (\kappa -2 ) \, \gs .
\ee
\eprop
\bp
Equation \eqref{eqhh4} can be written as
\be
\label{eqhh7}
\big(\log | \g (\ei)|^2 \big) \g (\ei)
\ \sim \ \kappa \g (\ei).
\ee
Differentiate both sides with respect to $\theta$. Writing $\g'$ for the derivative with respect to $z$, 
then \eqref{eqhh7} becomes
\[
\frac{iz \g'}{\g} \g \ - \ \frac{i \bar z \bar \g'}{\bar \g} \g + \log | \g |^2 (i z \g') \ \sim \ \kappa (i z \g') .
\]
So
\beq
 \log | \g |^2 (i z \g') &\ \sim \ & \kappa (i z \g') - i z \g' + i z^{n-1} \bar \g' \\
 & = &
 \kappa (i z \g') - i z \g'  + in \g - i z \g'. 
 \eeq
Therefore
\beq
\log | \g|^2 ( i \g - \frac{2i}{n} z \g') &\ \sim \ &
i \kappa \g - \frac{2i }{n} \left( \kappa ( z \g') -  z \g'  + n \g -  z \g' \right)
\\
 &\ = \ & i (\kappa -2) \left( \g - \frac{2}{n} z \g' \right) .
\eeq
This yields \eqref{eqhh6}.
\ep
It is plausible that the only polynomial satisfying \eqref{eqhh2} and \eqref{eqhh3} 
and with positive $0^{\rm th}$ and $n^{\rm th}$ coefficients is \eqref{eqhh5}, but
we cannot resolve whether this is true.

\section{A special case of the Entropy Conjecture}
\label{seci}
Self-inversive polynomials $p$ in $\pn$ can be written as
\be
\label{eqi1}
p \= q + \tq,
\ee
where $q$ is a polynomial in ${\mathcal P}_m$, with $m = \left \lfloor{\frac{n}{2}} \right \rfloor$.
Specifically, if $n$ is odd, then $q$ is an arbitrary polynomial in ${\mathcal P}_{\frac{n-1}{2}}$,
and defined by $$
q \= \Pi_{\frac{n-1}{2}}\  p ;
$$
if $n$ is even, then $q$ is a polynomial in ${\mathcal P}_{\frac{n}{2}}$ whose $\left( \frac{n}{2}\right)^{\rm th}$
coefficient is real (and half of the coefficient for $p$).
\bt
\label{thmi1}
Let $p$ be a self-inversive non-constant polynomial of degree $n$, and write $p$ as
in \eqref{eqi1}. If $q$ has no zeroes in the closed unit disk, then
Inequality \ref{eqa3} holds, with strict inequality unless $p$ is given by \eqref{eqa4}.
\et
\bp
Let us decompose the integral into two pieces, $I$ + $II$:
\be
\label{eqi2}
\int |p|^2 \log |p|^2 \= \int |p|^2 \log |q|^2 \ + \
\int |p|^2 \log | 1 + \tq/q|^2 .
\ee
To estimate  $I$, the first term on the right-hand side of \eqref{eqi2},
write 
\[
\int |p|^2 \log |q|^2 \= \int |q|^2 \log |q|^2 | | 1 + \tq/q|^2 .
\]
Note that 
\[
\int | 1 + \tq/q|^2 
\= \int 2 + 2 \Re \frac{\tq}{q} \= 2 ,
\]
since $\tq(0) = 0$.
So if we apply Jensen's inequality to the convex function $\Phi(x) = x \log x$
and the probability measure $\frac{1}{2}  | 1 + \tq/q|^2$, we get
\[
\int \Phi(|q|^2) \frac{1}{2}  | 1 + \tq/q|^2 \ \geq
\ \Phi \left( \int |q|^2  \frac{1}{2}  | 1 + \tq/q|^2 \right) .
\]
This gives
\[
\frac{1}{2} 
 \int |q|^2 \log |q|^2 | | 1 + \tq/q|^2
\ \geq \
\Phi(\frac{1}{2}) \= -\frac{1}{2} \log 2.
\]
Therefore we have
$I \geq - \log 2$.

To estimate $II$, first assume that $n$ is odd. Note that by the maximum principle, $\frac{\tq}{q}$ has modulus
less than one in the unit disk, so $\log ( 1 + \tq/q)$ is analytic on the unit disk and has only logarithmic singularities
on the unit circle, and therefore its Fourier series agrees with its Maclaurin series. Therefore
\[
\log ( 1 + \tq/q) \= \tq/q  + O(z^{n+1}) ,
\]
so
\beq
\int   |p|^2 \log |q|^2 &\= &
2 \Re \int |q + \tq |^2 \tq/q
\\
&=&  
2 \Re  \int ( 2 |q|^2 + z^n \bar q^2 + \bar z^n q^2 )(z^n \bar q/q)
\\
&=&
2 \Re \int 2 \bar q^2 z^n + z^{2n} \bar q^3 /q + |q|^2
\\ &=&
2 \Re \int |q|^2 \\
&=&
1.
\eeq

Now assume $n = 2m $ is even. Write
\[
q(z) \= a_0 + \dots + a_m z^m , \]
so
\[
\tq(z) \= a_m z^m + \bar a_{m-1} z^{m+1} + \dots + \bar a_0 z^n .
\]
When expanding $\log ( 1 + \tq/q)$
we get
\[
\log ( 1 + \tq/q) \= \frac{\tq}{q}  - \frac{1}{2} \frac{a_m^2}{a_0^2} z^n + O(z^{n+1}) .
\]
Therefore
\beq
\int   |p|^2 \log |q|^2 &\= &
2 \Re \int 2 \bar q^2 z^n + z^{2n} \bar q^3 /q + |q|^2 - \frac{1}{2}a_m^2 \\
&=&
3 a_m^2 + 1 \\
& \geq & 1.
\eeq
Therefore
\[
I + II \ \geq \ 1  - \log 2 ,
\]
as required.

Finally, note that the inequality for $I$ using Jensen's inequality is strict unless $|q|$ is constant.
\ep
 Note that a simple continuity argument applied to $q(rz)$ shows that \eqref{eqa3} holds provided
 $q$ has no zeroes in the open unit disk.
 
 \section{Baernstein's conjecture implies the entropy conjecture}
 \label{secj}
 
 Assume Baernstein's conjecture \ref{conb} holds. Let $Q(z) = (1+z^n)/\sqrt{2}$ and
 $p$ be any non-constant polynomial with all its roots on the unit circle, and with $\| p \|_2 = 1$.
 Let $s \leq 2$, and let $t = \frac{s}{s-1}$ be the conjugate exponent.
 Then, taking logarithms of  \ref{eqb2} we have
\[
\log \| p \|_s - \log \| p \|_t \ \leq \  
\log \| Q \|_s - \log \| Q \|_t .
\]
This means 
\be
\label{eqj1}
\frac{1}{s} \log \int |Q|^s - \frac{s-1}{s} \log \int |Q|^{s/(s-1)}
- \frac{1}{s} \log \int |p|^s + \frac{s-1}{s} \log \int |p|^{s/(s-1)} \geq 0.
 \ee
 Let $\Psi(s)$ denote the left-hand side of \eqref{eqj1}. We have that $\Psi(s) \geq 0$ for
 $s \leq 2$, and $\Psi(2) = 0$. Therefore $\Psi'(2) \leq 0$. 
 Calculating, we
 get 
 \beq
 \Psi'(s) &\=& 
 -\frac{1}{s^2} \log \int |Q|^s + \frac{1}{s} \frac{1}{\int |Q|^s} \int |Q|^s \log |Q|
 -\frac{1}{s^2}  \log \int |Q|^{\frac{s}{s-1}} \\
 &&
  +\frac{1}{s(s-1) } \frac{1}{ \int |Q|^{\frac{s}{(s-1}}}
 \int |Q|^{\frac{s}{s-1}} \log |Q| \\
 && +\ \frac{1}{s^2} \log \int |p|^s - \frac{1}{s} \frac{1}{\int |p|^s} \int |p|^s \log |p|
 + \frac{1}{s^2}  \log \int |p|^{\frac{s}{s-1}} \\
 &&
  -\ \frac{1}{s(s-1) } \frac{1}{ \int |p|^{\frac{s}{(s-1}}}
 \int |p|^{\frac{s}{s-1}} \log |p|
 \eeq
 Since both $p$ and $Q$ have $2$-norm $1$, we get that
 \[
 \Psi'(2) \= \int |Q|^2 \log |Q| - \int |p|^2 \log |p| .
 \]
 Since $\Psi'(2) \leq 0$, we get
 \[
 \int |p^2 \log |p|^2 \ \geq\  \int |Q|^2 \log |Q|^2 \= 1 - \log 2 ,
 \]
 which is \eqref{eqa3}.
 \ep

\bibliography{references}

\end{document}